\documentclass[11pt]{amsart}
\usepackage{times,epsf,amssymb,amsmath,hyperref, rlepsf}

\def\wh{\widehat}

\DeclareMathOperator{\spt}{spt}

\begin{document}

\newtheorem{thm}{Theorem}[section]
\newtheorem{lem}[thm]{Lemma}
\newtheorem{cor}[thm]{Corollary}

\theoremstyle{definition}
\newtheorem{defn}{Definition}[section]

\theoremstyle{remark}
\newtheorem{rmk}{Remark}[section]

\def\square{\hfill${\vcenter{\vbox{\hrule height.4pt \hbox{\vrule
width.4pt height7pt \kern7pt \vrule width.4pt} \hrule height.4pt}}}$}

\def\T{\mathcal T}

\newenvironment{pf}{{\it Proof:}\quad}{\square \vskip 12pt}

\title{Examples of Area Minimizing Surfaces in $3$-manifolds}
\author{Baris Coskunuzer}
\address{Koc University \\ Department of Mathematics \\ Sariyer, Istanbul 34450 Turkey}
\email{bcoskunuzer@ku.edu.tr}
\thanks{The author is partially supported by EU-FP7 Grant IRG-226062, TUBITAK Grant 109T685 and TUBA-GEBIP Award.}

\maketitle

%% User definitions:

\newcommand{\BR}{\mathbf{R}}
\newcommand{\BC}{\mathbf{C}}
\newcommand{\BZ}{\mathbf{Z}}
\newcommand{\BN}{\mathbf{N}}
\newcommand{\A}{\mathcal{A}}
\newcommand{\B}{\mathbf{B}}
\newcommand{\U}{\mathcal{U}}
\newcommand{\M}{\mathbf{M}}
\newcommand{\e}{\epsilon}
\newcommand{\de}{\delta}

\begin{abstract}

In this paper, we give some examples of area minimizing surfaces to clarify some well-known features of these surfaces in more general settings. The
first example is about Meeks-Yau's result on embeddedness of solution to the Plateau problem. We construct an example of simple closed curve in
$\BR^3$ which lies in the boundary of a mean convex domain in $\BR^3$, but the area minimizing disk in $\BR^3$ bounding this curve is not embedded.
Our second example shows that Brian White's boundary decomposition theorem does not extend when the ambient space have nontrivial homology. Our last
examples show that there are properly embedded absolutely area minimizing surfaces in a mean convex $3$-manifold $M$ such that while their boundaries
are disjoint, they intersect each other nontrivially.
\end{abstract}

\section{Introduction}

In this paper, we give examples of some area minimizing surfaces related with some well-known theorems. The first example is on Meeks-Yau's result on
embeddedness of solutions of the Plateau problem. The Plateau problem asks the existence of an area minimizing disk for a given curve in the ambient
manifold $M$. Meeks-Yau's famous result (Theorem \ref{Meeks-Yau}) says that for any simple closed curve in the boundary of a mean convex $3$-manifold
$M$, the solution to the Plateau problem \underline{in $M$} must be embedded. Since any convex body $C$ in $\BR^3$ is mean convex, and any solution
to the Plateau problem for a simple closed curve in $C$ must belong to $C$ because of the convexity, the result automatically implies that for any
extreme curve in $\BR^3$, the solution to the Plateau problem is embedded. However, our example shows that this is not the case for an $H$-extreme
curve $\Gamma$ in $\BR^3$, i.e. $\Gamma \subset \partial \Omega$ where $\Omega$ is mean convex in $\BR^3$. We construct an $H$-extreme curve $\Gamma$
in $\BR^3$ where the area minimizing disk $\Sigma$ in $\BR^3$ with $\partial \Sigma=\Gamma$ is not embedded (See Figure \ref{nonembedded2}).

Our second example shows that White's decomposition theorem for the absolutely area minimizing surfaces bounding curves with multiplicity does not
generalize to the manifolds with nontrivial second homology. It follows from \cite{Wh1} that if $\Gamma$ is a simple closed curve in $\BR^3$, and $T$
is an absolutely area minimizing surface with $\partial T = k\Gamma$ where $k$ is an integer greater than $1$, then $T = \sum_1^k T_i$ where each
$T_i$ is an absolutely area minimizing surface in $\BR^3$ with $\partial T_i = \Gamma$ (Theorem \ref{White}). This theorem naturally extends to
higher dimensions, and the orientable Riemannian manifolds with trivial second homology. Trivial homology plays crucial role in the proof, and we
give an example which shows that the decomposition theorem does not generalize to the $3$-manifolds with nontrivial second homology (See Figure
\ref{decomposition}).

Finally, in the last examples, we address the issue of intersections of absolutely area minimizing surfaces in mean convex $3$-manifolds. It is known
that if $\Gamma_1$ and $\Gamma_2$ are two disjoint simple closed curves in the boundary of a mean convex $3$-manifold $M$, then the area minimizing
disks they bound in $M$ must be disjoint, too \cite{MY2}. Our examples shows that the same statement is not true for absolutely area minimizing
surfaces. In particular, it is known that if $\Sigma_1$ and $\Sigma_2$ are two absolutely area minimizing surfaces with disjoint boundaries, then
they must also be disjoint, provided that $\Sigma_1$ and $\Sigma_2$ are homologous (rel. $\partial M$) \cite{Co}, \cite{Ha}. However, in the case
they are not homologous, and we construct absolutely area minimizing surfaces $\Sigma_1$ and $\Sigma_2$ in a mean convex $3$-manifold such that
$\partial \Sigma_i=\Gamma_i\subset \partial M$ and $\Gamma_1\cap\Gamma_2=\emptyset$, but $\Sigma_1\cap\Sigma_2\neq \emptyset$.

In particular, we can restate the result mentioned above (\cite{Co}, \cite{Ha}) in the following form: Let $M$ be a strictly mean convex
$3$-manifold, and let $\Sigma_1$ and $\Sigma_2$ be two absolutely area minimizing surfaces in $M$ with $\partial \Sigma_i = \Gamma_i \subset \partial
M$. Let $\Gamma_1\cap\Gamma_2 = \emptyset$, but $\Sigma_1 \cap \Sigma_2 \neq \emptyset$. Then, either $\Gamma_1$ is not homologous to $\Gamma_2$ in
$\partial M$, or $H_2(M)$ is not trivial. Our examples show that both situations are possible (See Figures \ref{intersecting1} and
\ref{intersecting2}).

The organization of the paper is as follows: In the next section, we give basic definitions and results which will be used in the following sections.
In Section 3, we describe the example about Meeks-Yau's embeddedness result. In Section 4, we give the example on White's decomposition theorem. In
the last section, we construct intersecting absolutely area minimizing surfaces with disjoint boundary.

\subsection{Acknowledgements:}
We would like to thank Brian White, Tolga Etgu and  Theodora Bourni for very useful remarks and conversations.

\section{Preliminaries}

In this section, we will give the basic definitions for the following sections.

\begin{defn} An {\em area minimizing disk} is a disk which has the smallest area among the disks with the same boundary.
An {\em absolutely area minimizing surface} is a surface which has the smallest area among all orientable surfaces (with no topological restriction)
with the same boundary.
\end{defn}

\begin{defn} \label{meanconvex} Let $M$ be a compact Riemannian $3$-manifold with boundary. Then $M$ is {\em mean convex} (or sufficiently convex) if the following conditions hold.

\begin{itemize}

\item $\partial M$ is piecewise smooth.

\item Each smooth subsurface of $\partial M$ has nonnegative curvature with respect to inward normal.

\item There exists a Riemannian manifold $N$ such that $M$ is isometric to a submanifold of $N$ and
each smooth subsurface $S$ of $\partial M$  extends to a smooth embedded surface $S'$ in $N$ such that $S' \cap M = S$.

\end{itemize}

\end{defn}

\begin{defn} A simple closed curve is an {\em extreme curve} if it is on the boundary of its convex hull.
A simple closed curve is called as {\em $H$-extreme curve} if it is a curve in the boundary of a mean convex manifold $M$.
\end{defn}

\section{Example I: Meeks-Yau Embeddedness Result}

The Plateau problem asks the existence of an area minimizing disk bounding a given curve in $\BR^3$. In other words, for a given simple closed curve
$\Gamma$ in $\BR^3$, does there exist a disk $\Sigma$ with the smallest area among the disks with boundary $\Gamma$, i.e. for any disk $D\subset
\BR^3$ with $\partial D = \Gamma$, $|D|\geq |\Sigma|$. This problem was solved by Douglas \cite{Do}, and Rado \cite{Ra} in early 1930s. Later, it was
generalized to homogeneously regular $3$-manifolds by Morrey \cite{Mo}. Then, the regularity (nonexistence of branch points) of these solutions was
shown by Osserman \cite{Os}, Gulliver \cite{Gu} and Alt \cite{Al}.

In the following decades, another version of the Plateau problem was investigated: Without any restriction on the genus of surface, does there exist
a smallest area surface bounding a given curve in $\BR^3$? In 1960s, the geometric measure theory techniques proved to be quite powerful, and De
Georgi, Federer-Fleming showed that for any simple closed curve $\Gamma$ in $\BR^3$ there exists an absolute area minimizing surface $S$ which
minimizes area among the surfaces (no topological restriction) with boundary $\Gamma$ \cite{Fe}. Moreover, any such surface is embedded in $\BR^3$.

When we go back to the original Plateau problem, in the 1970s, the question of embeddedness of the area minimizing disk was studied: For which curves
in $\BR^3$, the area minimizing disks are embedded? It was conjectured that if the simple closed curve in $\BR^3$ is extreme (lies in the boundary of
its convex hull), then the area minimizing disk spanning the curve must be embedded. After several partial results, Meeks and Yau proved the
conjecture: Any solution to the Plateau problem for an extreme curve must be embedded \cite{MY1}. Indeed, they proved more:

\begin{thm}\cite{MY2}, \cite{MY3} \label{Meeks-Yau}
Let $M$ be a compact, mean convex $3$-manifold, and $\Gamma\subset\partial M$ be a nullhomotopic simple closed curve. Then, there exists an area
minimizing disk $D\subset M$ with $\partial  D = \Gamma$. Moreover, all such disks are properly embedded in $M$.
\end{thm}

This theorem automatically implies the conjecture for extreme curves in $\BR^3$. If $\Gamma$ is an extreme curve in $\BR^3$, then let $M$ be the
convex hull of $\Gamma$, i.e. $M=CH(\Gamma)$. Since $CH(\Gamma)$ is a convex set in $\BR^3$, it automatically satisfies the mean convexity condition.
The theorem says that there is an area minimizing disk $\Sigma$ in $M$ with $\partial \Sigma = \Gamma$. Now, since $M$ is not just mean convex, but
also convex in $\BR^3$, the area minimizing disk $\Sigma$ in $M$ is \underline{also area minimizing in $\BR^3$}. In other words, $\Sigma$ has the
smallest area among the disks with boundary $\Gamma$ not only in $M$, but also in $\BR^3$. This shows that any solution to the Plateau problem in
$\BR^3$ must belong to $M$, hence by the theorem above, it must be embedded.

In this theorem, there is a subtle point. The theorem says that for a given simple closed nullhomotopic curve $\Gamma$ in $\partial M$, the area
minimizing disk \underline{in $M$} is embedded. This does not say that if $M$ is a mean convex domain in $\BR^3$, and $\Gamma$ is a simple closed
curve in $\partial M$, the area minimizing disk $D$ in $\BR^3$ with $\partial D =\Gamma$ is embedded. In other words, the theorem gives the
embeddedness of the minimizer in $M$, not the minimizer in $\BR^3$. In this paper, we will construct an explicit example for this difference. We note
that Spadaro has recently constructed an example of a simple closed curve $\Gamma$ on the boundary of a mean convex domain $\Omega$ in $\BR^3$ where
the area minimizing disk in $\BR^3$ with boundary $\Gamma$ does not belong to $\Omega$ \cite{Sp}. While the minimizers in $\BR^3$ are still embedded
in his examples, here we construct examples where the minimizers in $\BR^3$ are not embedded.

Now, we describe the example, a simple closed curve $\Gamma$ in $\BR^3$ which lies in the boundary of a mean convex domain $M \subset \BR^3$, but the
area minimizing disk $\Sigma$ in $\BR^3$ with $\partial \Sigma = \Gamma$ is not embedded.

\begin{figure}[h]

\relabelbox  {\epsfxsize=3.5in

\centerline{\epsfbox{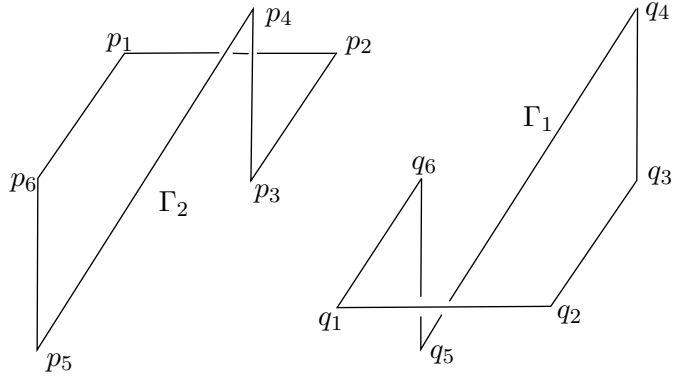}}}

\relabel{1}{$\Gamma_2$}

\relabel{2}{$\Gamma_1$}

\relabel{3}{$p_1$}

\relabel{4}{$p_2$}

\relabel{5}{$p_3$}

\relabel{6}{$p_4$}

\relabel{7}{$p_5$}

\relabel{8}{$p_6$}

\relabel{9}{$q_1$}

\relabel{10}{$q_2$}

\relabel{11}{$q_3$}

\relabel{12}{$q_4$}

\relabel{13}{$q_5$}

\relabel{14}{$q_6$}

\endrelabelbox

\caption{\label{nonembedded1} \small {The extreme curves $\Gamma_1$ and $\Gamma_2$.}}

\end{figure}

Consider $\BR^3$ with $xyz$ coordinate system. Define $\Gamma_1$ as follows:

Let $\alpha_1$ be the line segment in $xy$-plane connecting the points $p_1=(1,-1,0)$ and $p_2=(1,1,0)$. Let $\alpha_2$ be the line segment in
$xy$-plane connecting the points $p_2$ and $p_3=(\epsilon,1,0)$. Let $\alpha_3$ be the line segment in $\{x=\epsilon\}$-plane connecting the points
$p_3$ and $p_4=(\epsilon,1,C)$. Let $\alpha_4$ be the line segment in $\{x=\epsilon\}$-plane connecting the points $p_4$ and $p_5=(\epsilon,-1,-C)$.
Let $\alpha_5$ be the line segment in $\{x=\epsilon\}$-plane connecting the points $p_5$ and $p_6=(\epsilon,-1,0)$. Let $\alpha_6$ be the line
segment in $xy$-plane connecting the points $p_6$ and $p_1$. Then, let $\Gamma_1 = \bigcup_{i=1}^6 \alpha_i$ (See Figure \ref{nonembedded1}).

Similarly, define $\beta_i$ and $q_i$ as the reflection of $\alpha_i$ and $p_i$ with respect to $yz$ plane, and define $\Gamma_2 = \bigcup_{i=1}^6
\beta_i$, i.e. $\Gamma_2$ is the reflection of $\Gamma_1$ with respect to $yz$-plane.

Now, $\Gamma_1$ is an extreme curve as it is in boundary of a convex box $B_1$ with corners: $(1,-1,C), (1,1,C), p_4, (\epsilon, -1, C), (1,-1,-C),
(1,1,-C), (\epsilon, 1, -C), p_5$. Similarly, $\Gamma_2$ is in the boundary of the box $B_2$, the reflection of $B_1$ with respect to $yz$-plane, and
hence $\Gamma_2$ is an extreme curve, too. Without loss of generality, we can smooth out the corners of $\Gamma_1$ and $\Gamma_2$ so that $\Gamma_1$
and $\Gamma_2$ are smooth curves with $\Gamma_i\subset\partial B_i$.

Let $\Sigma_1$ and $\Sigma_2$ be the area minimizing disks in $\BR^3$ with $\partial \Sigma_i = \Gamma_i$. Since the boxes are convex,
$\Sigma_i\subset B_i$. Let $\gamma$ be the line segment between the points $p_5$ and $q_5$. Now, recall the bridge principle for stable minimal
surfaces.

\begin{lem} \cite[Theorem 7]{MY3} (Bridge Principle) Let $\Sigma_1$ and $\Sigma_2$ be two stable orientable compact minimally immersed surfaces in $\BR^n$
and $\gamma$ be a Jordan curve joining $\partial \Sigma_1$ and $\partial \Sigma_2$. Then for any tubular neighborhood of $\gamma$, we can find a
bridge pair joining $\partial \Sigma_1$ and $\partial \Sigma_2$ such that the new configuration is the boundary of a compact minimal surface which is
close to the union of $\Sigma_1$ and $\Sigma_2$ joint by a strip in the tubular neighborhood of $\gamma$.
\end{lem}

Then by applying the lemma above to $\Sigma_1$, $\Sigma_2$ and $\gamma$, we get a stable minimal surface $\widehat{\Sigma} \sim \Sigma_1\sharp_\gamma
\Sigma_2$ with the boundary $\widehat{\Gamma} = \Gamma_1\sharp_\gamma \Gamma_2$. In other words, $\widehat{\Sigma}$ is the stable minimal surface
close to the union of $\Sigma_1$ and $\Sigma_2$ joint by a strip near $\gamma$, and $\widehat{\Gamma}$ is its boundary.

Now, $\widehat{\Sigma}$ is an embedded stable minimal disk with smooth boundary $\widehat{\Gamma}$. By \cite{MY3} (Corollary 1 at page 159), there is
a mean convex neighborhood of $\widehat{\Sigma}$, say $M$, in $\BR^3$ such that $\partial \widehat{\Sigma} = \widehat{\Gamma} \subset \partial M$.
Hence, $\widehat{\Gamma}$ is a simple closed curve in $\BR^3$ which lies in the boundary of a mean convex domain $M$. Now, we claim that the area
minimizing disk $D$ in $\BR^3$ with $\partial D = \widehat{\Gamma}$ is not embedded.

Let $\tau$ be the square in the $xy$-plane with corners $p_1, p_2, q_1, q_2$. Let $E$ be the disk in $xy$-plane with boundary $\tau$. Clearly, $E$ is
the unique area minimizing disk with boundary $\tau$. Here, uniqueness comes from the foliation of $\BR^3$ by planes $\{z=t \ | \ t\in \BR \}$. If
there was another such area minimizing disk, it would have a tangential intersection with one of the planes at the maximum (or minimum) height, which
contradicts to the maximum principle.

First, note that by taking $C$ accordingly, we can make $|\widehat{\Sigma}| \sim |\Sigma_1|+|\Sigma_2|\sim 2(C+\sqrt{C^2+1})$ as large as we want,
where $|.|$ represents the area. Hence, by choosing $C$ sufficiently large, we can assume $|\widehat{\Sigma}| \gg |E|=4$.

Now, by modifying $E$ slightly, we will get another disk $\widehat{E}$ in $\BR^3$ with $\partial \widehat{E} = \widehat{\Gamma}$. Let $\eta_{1a}$ be
the line segment in the $yz$-plane connecting the points $(0,1,0)$ and $(0, 1, C)$. Let $\eta_{1b}$ be the line segment in the $yz$-plane connecting
the points $(0,1,C)$ and $(0,-1, -C)$. Let $\eta_1 = \eta_{1a}\cup \eta_{1b}$. Let $\eta_2$ be the line segment in the $yz$-plane connecting the
points $(0, -1, -C)$ and $(0,-1,0)$.

\begin{figure}[t]
\begin{center}
$\begin{array}{c@{\hspace{.2in}}c}

\relabelbox  {\epsfxsize=2.5in \epsfbox{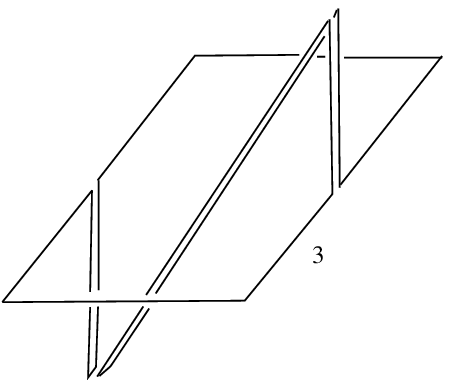}} \relabel{3}{$\Gamma^\epsilon$} \endrelabelbox &

\relabelbox  {\epsfxsize=2.5in \epsfbox{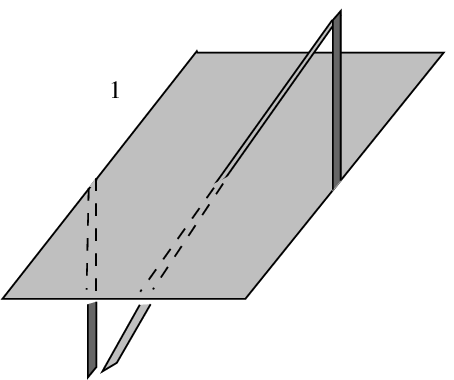}} \relabel{1}{$\wh{E}^\epsilon$} \endrelabelbox \\ [0.4cm]
\end{array}$

\end{center}
\caption{\label{nonembedded2}The H-extreme curve $\Gamma^\epsilon$ bounds a non-embedded area minimizing disk $\wh{E}^\epsilon$.}
\end{figure}

Let $S_{1a} = \{ (x,1,z) \ | \ x\in [-\epsilon, \epsilon],  \ z\in [0,C] \ \}$ be the strip near $\eta_{1a}$. Let $S_{1b} = \{ (x,y,Cy) \ | \ x\in
[-\epsilon, \epsilon],    y\in [-1,1] \}$ be the strip near $\eta_{1b}$. Finally, let $S_2 = \{ (x,-1,z) \ | \ x\in [-\epsilon, \epsilon], z\in
[-1,0] \}$. In other words, $S_{1a}$ is strip near $\eta_{1a}$ in $\{y=1\}$ plane, $S_{1b}$ is strip near $\eta_{1a}$ in $\{z=Cy\}$ plane, and $S_2$
is the strip in $\{y=-1\}$ plane with thicknesses $2\epsilon$. Let $S_1 = S_{1a} \cup S_{1b}$.

To get $\widehat{E}$, modify/trim the tips of the strips $S_1$ and $S_2$, say $\widehat{S}_1$ and $\widehat{S}_2$, so that $\partial \widehat{E} =
\widehat{\Gamma}$ where $\widehat{E} = E\cup \widehat{S}_1 \cup \widehat{S}_2$ (See Figure \ref{nonembedded2}). Clearly, $|\widehat{E}|\sim
|E|+|S_1|+|S_2|= 4+2\epsilon(C+\sqrt{C^2+1})$. Hence, by fixing sufficiently large $C_0$, there exists $\epsilon_0>0$ such that
$|\widehat{\Sigma}^\epsilon|>|\widehat{E}^\epsilon|$ for any $\epsilon<\epsilon_0$. This shows that the area minimizing disk in $\BR^3$ bounding
$\widehat{\Gamma}^\epsilon$ does not lie in $M$.

From now on, we fix a sufficiently large $C_0$ as above, and $\widehat{\Sigma}^\epsilon,\widehat{E}^\epsilon$ and $\widehat{\Gamma}_\epsilon$
represents the corresponding disks and curve for $\epsilon<\epsilon_0$. Hence, $\partial \widehat{\Sigma}^\epsilon = \partial \widehat{E}^\epsilon =
\widehat{\Gamma}_\epsilon$, and  $|\widehat{\Sigma}^\epsilon| \gg |\widehat{E}^\epsilon|$ for any $\epsilon<\epsilon_0$.\\

\noindent {\bf Claim:} There is a sufficiently small $\epsilon<\epsilon_0$  such that the area minimizing disk $D^\epsilon$ in $\BR^3$ with $\partial
D^\epsilon= \widehat{\Gamma}^\epsilon$ is not embedded.\\

\begin{pf} First, we will show that $D^\epsilon$ is very close to $\widehat{E}^\epsilon$ for small $\epsilon>0$. In particular, we claim that for any
$\rho>0$, there exists $\epsilon<\epsilon_0$ such that $d(D^\epsilon,\widehat{E}^\epsilon)<\rho$ where the distance $d$ represents the Hausdorff
distance, and $D^\epsilon$ is the area minimizing disk in $\BR^3$ with $\partial D^\epsilon = \widehat{\Gamma}^\epsilon$.

Assume on the contrary that there is a $\rho_0>0$ such that $d(D^\epsilon,\widehat{E}^\epsilon)\geq \rho_0$ for any $\epsilon>0$. Let
$\epsilon_i\searrow 0$ be a sequence converging to $0$. Let $T_i=[\![D^{\e_i}]\!]$ be the corresponding currents. Then, $\partial T_i = Y_i=
[\![\widehat{\Gamma}^{\e_i}]\!]$. Let $T_0=[\![E]\!]$ and $Y_0=[\![\tau]\!]$.

Let $\M$ represent the mass of a current. Then, $\M(T_i)=|D^{\e_i}|\leq |\widehat{E}^{\e_i}| = 4+2\epsilon(C_0+\sqrt{C_0^2+1})$ for any $i$ as
$D^{\e_i}$ is the area minimizing disk, and $\widehat{E}^{\e_i}$ is a disk with the same boundary $\widehat{\Gamma}^{\e_i}$. Hence, the currents
$T_i$ have uniformly bounded masses, and thus by the Federer-Fleming compactness theorem \cite{Fe}, after passing to a subsequence, $T_i\to T$, in
the sense of currents, where $T$ is an integral current.

We claim that $T=T_0=[\![E]\!]$. As $M(T_i)\to 4$, we get $M(T)=4$. By construction $Y_i\to Y_0$ in the sense of currents, and hence $\partial T =
Y_0$. Note that $E$ is the unique absolutely area minimizing surface in $\BR^3$ with $\partial E = \tau$. As $M(T)=|E|$ with $\partial T = Y_0$, this
implies $T=T_0$. This proves that $T_i\to T_0$ in the sense of currents.

Recall that by assumption, there is $\rho_0>0$ such that $d(D^{\e_i},\widehat{E}^{\e_i})\geq \rho_0$ for any $i$. In other words, there exists
$x_i\in D^{\e_i}$ such that $d(x_i, \widehat{E}^{\e_i})\geq \rho_0$ for any $i$. After passing to a subsequence, we get $x_i\to x_0 \in \BR^3$, and
by construction $d(x_0, E \cup \eta_1\cup \eta_2 )\geq \rho_0$ as $E \cup \eta_1\cup \eta_2 \subset \widehat{E}^{\e_i}$ for any $i$. Since $x_i\to
x_0$, there exists $i_0$, such that for any $ i\ge i_0$ we have $\B_{\rho_0/4}(x_i)\subset \B_{\rho_0/2}(x_0)\subset \B_{\rho_0}(x_i)$ where
$\B_\rho(x)$ denotes the ball of radius $\rho$ and centered at $x$ in $\BR^3$.

Hence, using the monotonicity formula \cite{Si} and the fact that $x_i\in\spt T_i\setminus \spt\partial T_i$, we get $\mu_{T_i}(\B_{\frac {\rho_0}
{2}}(x_0))\ge \mu_{T_i}(\B_{\frac {\rho_0} {4}}(x_i))\ge \pi\left(\frac{\rho_0}{4}\right)^2$

Now, by using the measure convergence, we have that $$ \mu_{T_0}(\B_{\frac {\rho_0} {2}}(x_0)) = \lim_i\mu_{T_i}(\B_{\frac {\rho_0} {2}}(x_0))\geq
\pi\left(\frac{\rho_0}{4}\right)^2$$ and since $\B_{\frac {\rho_0} {2}}(x_0)\subset \B_{\rho_0}(x_i)$ $$\mu_{T_0}(\B_{\rho_0}(x_i)) \geq
\mu_{T_0}(\B_{\frac {\rho_0} {2}}(x_0))\geq \pi\left(\frac{\rho_0}{4}\right)^2$$ which implies that $\spt T_0 \cap \B_{\rho_0}(x_i)\ne \emptyset$.
However, by assumption $d (x_i, E)\geq \rho_0$, and this is contradiction.

Hence, for any $\rho>0$, there is a sufficiently small $\epsilon>0$ such that the area minimizing disk $D^\epsilon$ is $\rho$-close to
$\widehat{E}^\epsilon$. Since $\widehat{E}^\epsilon$ has transversal self-intersection, this implies $D^\epsilon$ is not embedded for sufficiently
small $\epsilon$. This finishes the proof of the claim.
\end{pf}

Hence, the claim above implies the area minimizing disk $D^\e$ in $\BR^3$ which $\Gamma^\epsilon$ bounds is not embedded, even though
$\Gamma^\epsilon$ lies in the boundary of a mean convex domain $M^\epsilon$ in $\BR^3$. In other words, Meeks-Yau's result (Theorem \ref{Meeks-Yau})
says that for an $H$-extreme curve $\Gamma\subset \partial M$, the solution to the Plateau problem \underline{in the mean convex manifold $M$} is
embedded, and this implies the solution to the Plateau problem for an extreme curve in $\BR^3$ is embedded. However, it does not say that the
solution to the Plateau problem for an $H$-extreme curve \underline{in $\BR^3$} is embedded, and the example above shows that this is not true.

\begin{rmk} Let $\wh{\Sigma}^\e$ and $D^\e$ be the area minimizing disks in the example above with $\partial \wh{\Sigma}^\e = \partial D^\e = \wh{\Gamma}^\e$
($\wh{\Sigma}^\e$ is the minimizer in $M^\e$ while $D^\e$ is minimizer in $\BR^3$). Notice that by choosing $C$ and $\epsilon$ accordingly in the
construction, the ratio $\frac{|\wh{\Sigma}^\e|}{|D^\e|}$ between the the area minimizing disks $\wh{\Sigma}^\e$ and $D^\e$ can be made as large as
we want. In other words, while the area of the minimizer in the ambient space $\BR^3$ is very small, the minimizer in the mean convex manifold $M^\e
\subset \BR^3$ can have very large area.
\end{rmk}

\section{Example II: White's Decomposition Theorem}

White's decomposition theorem for absolutely area-minimizing hypersurfaces at boundaries with multiplicity states the following:

\begin{thm} \cite{Wh1} \label{White} Let $\Gamma$ be a codimension-$2$ smooth submanifold in $\BR^{n+1}$, and $T$ is an absolutely area minimizing hypersurface with $\partial T =
k\Gamma$ where $k>1$. Then, $T = \sum_1^k T_i$ where each $T_i$ is an absolutely area minimizing hypersurface with $\partial T_i = \Gamma$.
\end{thm}

In particular, this theorem extends naturally to the orientable Riemannian manifolds with trivial second homology. In this paper, by constructing an
explicit example, we will show that the decomposition theorem does not generalize to the orientable manifolds with nontrivial second homology. Note
that our example is $3$-dimensional, but by using similar arguments, it can be extended to higher dimensions.

\begin{figure}[h]

\relabelbox  {\epsfxsize=3.5in

\centerline{\epsfbox{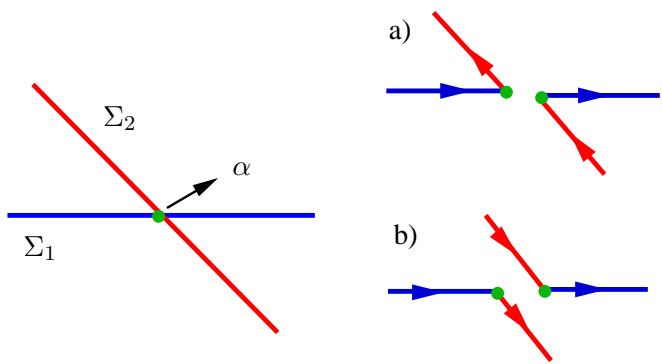}}}

\relabel{1}{$\Sigma_1$}

\relabel{2}{$\Sigma_2$}

\relabel{3}{$\alpha$}

\relabel{4}{a)}

\relabel{5}{b)}

\endrelabelbox

\caption{\label{surgery} \small One dimensional surgery case is pictured. $\Sigma_1$ and $\Sigma_2$ intersects along a simple closed curve $\alpha$.
After cutting both surfaces along $\alpha$, we can glue the surfaces as shown in the figure right. Depending on the orientations of $\Sigma_1$ and
$\Sigma_2$, surgery $a)$ or surgery $b)$ gives us an oriented surface. }

\end{figure}

An important observation about this theorem is the following: If this decomposition theorem was true in general setting, it would imply that the
absolutely area minimizing surfaces with the same boundary in such a manifold would be disjoint. This is because if $\Sigma_1$ and $\Sigma_2$ are two
absolutely area minimizing surfaces with $\partial \Sigma_i = \Gamma$, and $\Sigma_1\cap \Sigma_2 \neq \emptyset$, then we can make a surgery along
the intersection curve $\alpha$ (See Figure \ref{surgery}).

First, note that $\Sigma_1$ and $\Sigma_2$ cannot be separating each other. This is because if the two surfaces are separating each other, then we
get a contradiction by a swaping argument as follows: If $\Sigma_1$ and $\Sigma_2$ are separating each other, then $\Sigma_1 -\Sigma_2 = S_1^+\cup
S_1^-$ and $\Sigma_2 -\Sigma_1 = S_2^+\cup S_2^-$ where $S_i^+$ be the component in $\Sigma_i$ containing the boundary $\Gamma$, and $S_i^-$ be the
other components in $\Sigma_i$. Note that $|S_1^-|=|S_2^-|$ as they are also absolutely area minimizing surfaces with the same boundary. Then define
$\Sigma_1' = S_1^+ \cup S_2^-$ and $\Sigma_2'= S_2^+\cup S_1^-$. Notice that the new surfaces $\Sigma_1'$ and $\Sigma_2'$ have the same area with
$|\Sigma_1|=|\Sigma_2|$. However, new surfaces contain folding curves $\Sigma_1\cap\Sigma_2$. We get a contradiction as we can get a smaller area
surfaces by pushing the surfaces to the convex side along folding curve \cite{MY2}. Note that the same contradiction can be obtained by using the
regularity theorem for absolutely area minimizing surfaces as the absolutely area minimizing surfaces must be smooth \cite{Fe}. This implies that the
surfaces can not separate each other.

Note that if two surfaces were homologous, either the intersection would be empty, or the surfaces would separate each other. Since $\Sigma_1$ and
$\Sigma_2$ do not separate each other, by doing surgery (choosing the correct sides to match the orientation) along $\alpha$, we get a connected
oriented surface $\Sigma$ with the area $|\Sigma_1|+|\Sigma_2|$ and with $\partial \Sigma = 2\Gamma$ (See Figure \ref{surgery}). $\Sigma$ is an
absolutely area minimizing surface for the boundary $2\Gamma$ as the decomposition theorem implies that the smallest area surface with boundary
$2\Gamma$ has area $|\Sigma_1|+|\Sigma_2| = 2|\Sigma_i|$. This is a contradiction as before, since we have a folding curve along $\alpha$ in the
absolutely area minimizing surface $\Sigma$.

As we will see in the following section, this is not the case (See Remark \ref{samebdry}). Hence, the decomposition theorem is not true for the
oriented $3$-manifolds with nontrivial second homology. We will construct an example of a \underline{connected}, oriented absolutely area minimizing
surface $\Sigma$ with $\partial \Sigma = 2\Gamma$ with less area than $|\Sigma_1|+|\Sigma_2|$, which shows that the decomposition theorem does not
extend to manifolds with nontrivial homology.

Now, we construct the example. First, to give the basic idea, we will start the topological construction where the area minimizing surface is not
smooth. Then, we will modify the original example to make the ambient space strictly mean convex and the surfaces are smooth (See Remark \ref{NMC}).

\begin{figure}[h]

\relabelbox  {\epsfxsize=3.5in

\centerline{\epsfbox{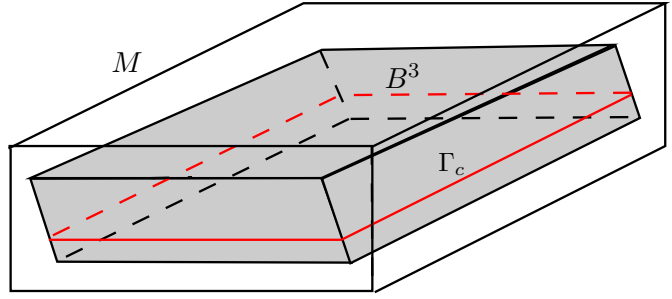}}}

\relabel{1}{$M$}

\relabel{2}{$B^3$}

\relabel{3}{\small $\Gamma_c$}

\endrelabelbox

\caption{$M$ is obtained by removing $B^3$ from the $3$-torus. $\Gamma_c$ is a simple closed curve in $\partial M = \partial B^3$.}

\end{figure}

Let $T^3$ be the $3$-torus obtained by identifying the opposite faces of the rectangular box of dimensions $[0,1]\times[0,1]\times [0,h]$. Take the
induced flat metric on $T^3$. Let $B^3$ be the $3$-ball in the shape of a parallelepiped in the $T^3$ with square base of dimensions
$[\delta,1-\delta]\times[\delta,1-\delta]$ and height $\frac{h}{3}$. We also assume that the parallelepiped has slope $\theta_0$ in $x$-direction
where $\tan{\theta_0}<1/6$. In particular, the corners of the parallelepiped is as follows: The vertices of the bottom square are $a_1=(\de,
\de,\frac{h}{3})$ , $a_2=(\de, 1-\de,\frac{h}{3})$ , $a_3=(1-\de, 1-\de,\frac{h}{3})$ , $a_4=(1-\de, \de,\frac{h}{3})$. Let $\sigma=
\frac{h}{3}\cot{\theta_0}>2h$. Then, by shifting all the bottom vertices $\sigma$ (assume $\sigma<\delta$) in the negative $x$-direction, we get the
vertices of the top square: $b_1=(\de-\sigma, \de,\frac{2h}{3})$ , $b_2=(\de-\sigma, 1-\de,\frac{2h}{3})$ , $b_3=(1-\de-\sigma, 1-\de,\frac{2h}{3})$
, $b_4=(1-\de-\sigma, \de,\frac{2h}{3})$ Now, define the ambient space $M$ as $T^3-B^3$ where $h$ and $\delta$ are to be determined later.

For $\frac{h}{3}<c<\frac{2h}{3}$, let $\Gamma_c$ be the simple closed curve in $\partial M$ with $\Gamma_c = \partial M \cap T_c$ where $T_c$ is the
$2$-torus corresponding to $\{z=c\}$ square in $T^3$ before identification. In other words, $\Gamma_c$ is the union of four line segments, i.e.
$\Gamma= \tau_1\cup\tau_2\cup\tau_3\cup\tau_4$ where $\tau_1=\{\delta\}\times[\delta,1-\delta]\times \{c\}$, $\tau_2= [\delta,1-\delta]\times
\{1-\delta\} \times \{c\}$, $\tau_3= \{1-\delta\}\times[\delta,1-\delta]\times \{c\}$, and $\tau_4=[\delta,1-\delta]\times \{\delta\} \times \{c\}$.

\begin{figure}[b]

\relabelbox  {\epsfxsize=3.5in

\centerline{\epsfbox{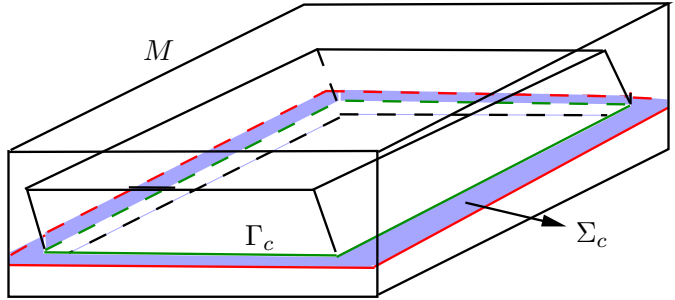}}}

\relabel{1}{$M$}

\relabel{2}{$\Gamma_c$}

\relabel{3}{$\Sigma_c$}

\endrelabelbox

\caption{\label{horizontal} $\Sigma_c$ is the absolutely area minimizing surface in $M$ with $\partial \Sigma_c = \Gamma_c$.}

\end{figure}

Let $\Sigma_c$ be the surface obtained by the intersecting the $2$-torus $T_c$ with $M$. Hence, $\Sigma_c$ is a $2$-torus where a disk ($T_c\cap
B^3$) is removed and $\partial \Sigma_c =\Gamma_c$ (See Figure \ref{horizontal}). Observe that either $\Sigma_c$ is the absolutely area minimizing
surface in $M$ with $\partial \Sigma_c = \Gamma_c$ or the absolutely area minimizing surface in $M$ with boundary $\Gamma_c$ completely lies in
$\partial M$. This is because the family of horizontal $2$-tori $\{T_c \}$ foliates the flat torus $T^3$ by minimal surfaces. If there was another
area minimizing surface $S$ in $M$ with $\partial S=\Gamma_c$, which does not completely lie in $\partial M$, $S \subset T^3$ must have a tangential
intersection with one of horizontal $2$-tori $T_{c'}$ by lying in one side, which contradicts to the maximum principle.

Assuming $c<\frac{h}{2}$, let $D_c$ be the smaller disk in $\partial M$ with $\partial D_c = \Gamma_c$. Then, the area of $D_c$ would be $x^2 +
4x\frac{c- h/3}{\sin{\theta_0}}$ where $x=1-2\delta$, the side of the square base of $\partial M$. The area of $\Sigma_c$ is $1-x^2$ for any $c\in
[\frac{h}{3},\frac{2h}{3}]$. Hence, when $x^2 + 4x\frac{c- h/3}{\sin{\theta_0}}> 1-x^2$, $\Sigma_c$ is the absolutely area minimizing surface, and
when $x^2 + 4x\frac{c- h/3}{\sin{\theta_0}} < 1-x^2$, $D_c$ is the absolutely area minimizing surface in $M$ with boundary $\Gamma_c$ by the
discussion above. Also, for $c_o\in [\frac{h}{3},\frac{2h}{3}]$ with $x^2 + 4x\frac{c_o- h/3}{\sin{\theta_0}} = 1-x^2$, both $\Sigma_{c_o}$ and
$D_{c_o}$ would be absolutely area minimizing surfaces with boundary $\Gamma_{c_o}$. Notice that by choosing $x=1-2\de$ sufficiently close to
$1/{\sqrt{2}}$ from below, we can choose $h$ as small as we want, and we would have a solution $c_o$ to the equation $x^2 + 4x\frac{c_o-
h/3}{\sin{\theta_0}} = 1-x^2$ in $[\frac{h}{3},\frac{2h}{3}]$. Now, we will show that there is a connected absolutely area minimizing surface
$\Sigma$ with $\partial \Sigma = 2\Gamma_{c_o}$.

\begin{figure}[b]

\relabelbox  {\epsfxsize=5in

\centerline{\epsfbox{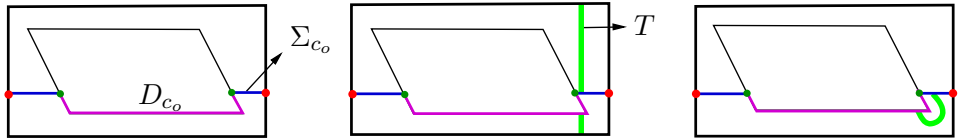}}}

\relabel{1}{$D_{c_o}$}

\relabel{2}{$\Sigma_{c_o}$}

\relabel{3}{$T$}

\endrelabelbox

\caption{\label{decomposition} \small The green dots represents $\Gamma_{c_o}$. The red dots represents the identification in the torus.
$\Sigma_{c_o}$ (blue surface) and $D_{c_o}$ (purple surface) are both absolutely area minimizing surfaces bounding $\Gamma_{c_o}$ in $M$. In the
middle figure, the tube $T$ connecting $D_{c_o}$ and $\Sigma_{c_o}$ gives the connected surface $\Sigma$ with $\partial \Sigma = 2\Gamma_{c_o}$ with
less area. In the figure right, we see that if we place the green tube in that way, the surface we are going to have will have boundary
$0.\Gamma_{c_o}=\emptyset$ as $D_{c_o}$ and $\Sigma_{c_o}$ must be oppositely oriented at the beginning to have an oriented surface after surgery in
this configuration. }

\end{figure}

Now, as in the figure, remove a disk $O_1$ from $\Sigma_{c_o}$ and remove a disk $O_2$ from $D_{c_o}$ where both $O_1$ and $O_2$ are disks with
radius $\epsilon$ where $h<\epsilon<\frac{\sigma}{2}=\frac{h}{6}\cot{\theta_0}$. In particular, let $O_1$ have center
$(1-\delta-\frac{\sigma}{2},\frac{1}{2}, c_o)$, and $O_2$ have center $(1-\delta-\frac{\sigma}{2},\frac{1}{2}, \frac{h}{3})$. Let $T$ be the cylinder
in $M$ with boundary $\partial O_1 \cup\partial O_2$ as in the Figure \ref{decomposition} (middle). Then, the area of the tube $T$ would be
$2.\pi.\epsilon.(h-(c_o-\frac{h}{3}))=2.\pi.\epsilon.(\frac{4h}{3}-c_o)$, while the area of the disks $O_1\cup O_2$ would be $2.\pi.\epsilon^2$.
Hence, for $\frac{4h}{3}-c_o< h<\epsilon$, $|T|<|O_1|+|O_2|$, and the new surface $\Sigma' = (\Sigma_{c_o} - O_1) \cup (D_{c_o} -O_2) \cup T$ would
have less area than the $\Sigma_{c_o}\cup D_{c_o}$. Moreover, by the choice of $T$, $\Sigma '$ is an orientable surface with $\partial \Sigma ' =
2\Gamma_{c_o}$ (See Remark \ref{handle}). Recall that if the decomposition theorem was true in this case, then the oriented absolutely area
minimizing surface $S$ with $\partial S = 2\Gamma_{c_o}$ would decompose as $ S= S_1+S_2$ where $S_i$ is an absolutely area minimizing surface with
$\partial S_i = \Gamma_{c_o}$. This would imply $|S| = |S_1|+|S_2|=2|S_1|$. However, there is an oriented surface $\Sigma '$ with $\partial \Sigma '
= 2\Gamma_{c_o}$ and $|\Sigma '| < |\Sigma_{c_o}| +|D_{c_o}|= 2|\Sigma_{c_o}|=2|D_{c_o}|$. This proves that the decomposition theorem is not valid
for ambient manifold $M$.

\begin{rmk} \label{handle} Note that here the choice of the handle, the cylinder $T$, in the construction is very important to get an oriented surface with boundary
$2\Gamma_{c_o}$. The other choice of the handle $T'$ as in the Figure \ref{decomposition} (right) would gave us another oriented surface with
boundary $0.\Gamma_{c_o}=\emptyset$ as we need to reverse the orientation on $D_{c_o}$ or $\Sigma_{c_o}$ to have oriented $T'$.
\end{rmk}

Now, in order to get smooth examples, we will modify the metric on $M$. To do this, we will change the metric $g$ on $M$ to $\wh{g}$ so that $M$ with
this new metric $\wh{g}$, say $\wh{M}$, will be strictly mean convex (See Remark \ref{NMC}). Then, the absolutely area minimizing surfaces in
$\wh{M}$ would be smoothly embedded (\cite{ASS}, \cite{Wh2}). Clearly, $M$ above with the induced flat metric is not mean convex as the dihedral
angles at the boundary are greater than $\pi$ (See \ref{meanconvex}, condition 3).

Note that as the mean convexity is a local condition, it will suffice to change the metric only near the boundary $\partial M$. Hence, the new metric
$\wh{g}$ will be same with $g$ everywhere on $M$ except a small neighborhood of $\partial M$, say $N_\xi(\partial M)$. First, change the metric in
$N_{\frac{\xi}{2}}(\partial M)$ so that it is isometric to $N_{\frac{\xi}{2}}(\partial B_{r_o})$ in $B_{r_o}$ where $B_{r_o}$ is the closed ball of
radius $r_o$ in $\BR^3$ with $|\partial B_{r_o}| = |\partial M|$. Here, we chose $r_o$ with $|\partial B_{r_o}| = |\partial M|$, since after
modification of the metric, we want the curve $\Gamma_{c}\subset \partial \wh{M}$ constructed above to have similar features as before. In
particular, in order to employ the construction above, after the modification of the metric, we would like to have a $c_o$ with for $c<c_o$, the
absolutely area minimizing surface bounding $\Gamma_c$ is the disk $D_{c_o}$ near boundary, and for $c>c_o$, the absolutely area minimizing surface
bounding $\Gamma_c$ is the punctured torus $\Sigma_c$.

After making the metric isometric to $N_{\frac{\xi}{2}}(\partial B_{r_o})$ in $N_{\frac{\xi}{2}}(\partial M)$, then to make the new metric smooth on
$M$, use partition of unity on the part $N_\xi(\partial M) - N_{\frac{\xi}{2}}(\partial M)$, which is homeomorphic to $S^2\times I$ so that it has
small cross-sectional area (cross-sectional annuli $\gamma\times I$ have small area $\sim \xi . |\gamma|$). We ask this to keep the areas of
$\Sigma_c$ small as before. Hence, we get a smooth metric on $M$ so that $M$ is strictly mean convex. Note that we can construct $\wh{g}$ as
rotationally symmetric in the $xy$-direction because of the setting.

We will call $M$ with the new metric $\wh{g}$ as $\widehat{M}$ for short. Note that we just changed the original metric in a very small neighborhood
of the boundary, $N_\xi(\partial M)$. We will use the same coordinates as before. Let $\{\Gamma_c\}$ be the family of simple closed curves in
$\partial \widehat{M}$ as above. Consider the absolutely area minimizing surfaces in $\widehat{M}$ with boundary $\Gamma_c$. By the construction of
the metric, these absolutely area minimizing surfaces would be close to $\Sigma_c$ or $D_c$ depending on $c$ as before. Call the absolutely area
minimizing surfaces in $\wh{M}$ close to $\Sigma_c$ and $D_c$ as $\widehat{\Sigma}_c$ and $\widehat{D}_c$ respectively. Notice that {\underline the
area} of $\wh{\Sigma}_c$ and $\wh{D}_c$ changes continuously with respect to $c$, as $\Gamma_i \to \Gamma$ the area of the annulus $A_i$ between
$\Gamma_i$ and $\Gamma$ goes to $0$, i.e. $|A_i|\to 0$. Hence, for some $c_o$, $|\wh{\Sigma}_{c_o}|=|\wh{D}_{c_o}|$. Now, we will imitate the
construction above, but this time we will get a smooth, oriented connected surface at the end.

Like before, remove a disk $\wh{O}_1$ from $\wh{\Sigma}_{c_o}$ and remove a disk $\wh{O}_2$ from $\wh{D}_{c_o}$ where $\wh{O}_1$ and $\wh{O}_2$ are
both disks with radius $\epsilon$. As before, we can choose the disks close in horizontal direction. Then, there is an area minimizing cylinder
$\wh{T}$ with the boundary $\partial \wh{O}_1 \cup \partial \wh{O}_2$ by \cite{MY2}, as $|\wh{T}| \sim |T|= 2.\pi.\epsilon.(h-c)$ and
$|\wh{O}_1|+|\wh{O}_2| \sim |O_1|+|O_2|= 2.\pi.\epsilon^2$. Let $\wh{\Sigma}' = (\wh{\Sigma}_{c_o}-\wh{O}_1) \cup (\wh{D}_{c_o} - \wh{O}_2) \cup
\wh{T}$. Then, like before $|\wh{\Sigma}'| < 2|\wh{D}_{c_o}| = |\wh{\Sigma}_{c_o}| + |\wh{D}_{c_o}|$ and $\partial \wh{\Sigma}' = 2\Gamma_{c_o}$.

\begin{rmk} If we took $\Sigma_1$ and $\Sigma_2$ in the example as homologous surfaces, the oriented surface we would get after the surgery would
have boundary $0.\Gamma$ not $2\Gamma$. This is because if $\Sigma_1$ and $\Sigma_2$ are homologous in $M$, then $\Sigma_1\cup\Sigma_2$ would
separate $M$ into $\Omega^+ \cup \Omega^-$. Hence, for any choice of the handle $T$ ($T\subset \Omega^+$ or $T\subset\Omega^-$), if $T$ is oriented
surface with $\partial T = \partial O_1 \cup \partial O_2$ where $O_i$ is a small disk in $\Sigma_i$ as above, then $T$ would be homologous to
$(\Sigma_1-O_1) +(\Sigma_2 -O_2)$. This would force $\Sigma_1$ and $\Sigma_2$ to be oppositely oriented, and we get an oriented connected surface
$\Sigma'=(\Sigma_1-O_1) \cup (\Sigma_2 -O_2) \cup T$ with less area. However, in this situation, $\partial \Sigma' = 0.\Gamma_{c_o} =\emptyset$.

Intuitively, when $\Sigma_1$ and $\Sigma_2$ are both oriented so that $\partial \Sigma_1=\partial\Sigma_2 = + \Gamma$, the "upper sides" of
$\Sigma_1$ and $\Sigma_2$ look into the same direction. To keep this orientations after surgery, the surgery tube $T$ must start from the "upper
side" of $\Sigma_1$ and end in the "upper side" of $\Sigma_2$. However, If $\Sigma_1 \cup \Sigma_2$ bounds the domain $\Omega^+$, and $T\subset
\Omega^+$, this is not possible as $T$ starts from the upper side of $\Sigma_1$ but ends in "bottom side" of $\Sigma_2$. Hence, to have $\Sigma'$
oriented, we must reverse the orientation on $\Sigma_2$, which forces $\partial \Sigma'=\Gamma - \Gamma = \emptyset$ (See Figure \ref{orientation}
(left)).

\begin{figure}[b]

\relabelbox  {\epsfxsize=5in

\centerline{\epsfbox{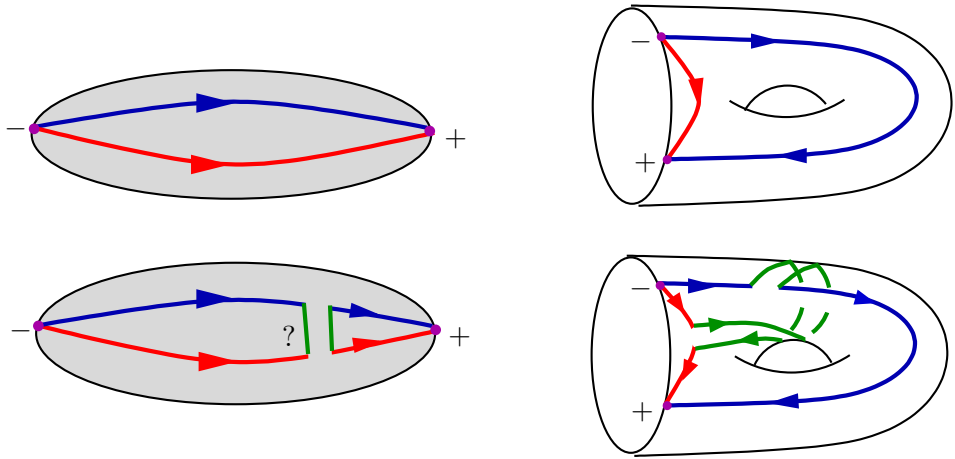}}}

\relabel{1}{$-$}

\relabel{2}{$+$}

\relabel{3}{$-$}

\relabel{4}{$+$}

\relabel{5}{$?$}

\relabel{6}{$-$}

\relabel{7}{$+$}

\relabel{8}{$-$}

\relabel{9}{$+$}

\endrelabelbox

\caption{\label{orientation} \small  If the surfaces $\Sigma_1$ and $\Sigma_2$ are homologous, we may not make a surgery (green "tube") consisting
with the original orientations given on $\Sigma_1$ and $\Sigma_2$ (left). However, if $\Sigma_1$ and $\Sigma_2$ are not homologous, there is a
surgery consistent with the original orientations for any given orientations on $\Sigma_1$ and $\Sigma_2$ so that $\partial \Sigma_1 \sharp \Sigma_2
= 2\Gamma$ (right). }

\end{figure}

However, if $\Sigma_1$ and $\Sigma_2$ are not homologous, then $\Sigma_1 \cup \Sigma_2$ does not separate $M$, and the tube $T$ is free to land in
the "upper side" $\Sigma_2$, which gives the surface $\Sigma'$ with the desired orientation (See Figure \ref{orientation} (right)).
\end{rmk}

\section{Example III: Intersections of Absolutely Area Minimizing Surfaces}

In this section, we will give explicit examples describing when the two absolutely area minimizing surfaces with disjoint boundaries can intersect.
Let $M$ be a mean convex manifold. Let $\Sigma_1$ and $\Sigma_2$ be two absolutely area minimizing surfaces in $M$ with $\partial \Sigma_i =
\Gamma_i$, where $\Gamma_1$ and $\Gamma_2$ are disjoint simple closed curves in $\partial M$. Then, it is known that if $\Sigma_1$ and $\Sigma_2$ are
homologous in $M$ (relative to $\partial M$), then they must be disjoint (\cite{Co}, Lemma 4.1), (\cite{Ha}, Theorem 2.3).

Hence, if two such absolutely area minimizing surfaces with disjoint boundaries intersect each other, they cannot be homologous. This gives us two
situations about the ambient space and boundary curves. The first situation is $H_2(M)$ is not trivial, and $\Sigma_1$ and $\Sigma_2$ belongs to
different homology classes in $H_2(M,\partial M)$. The second situation is $H_2(M)$ is trivial, but $\Gamma_1$ and $\Gamma_2$ are not homologous in
$\partial M$. In this paper, we will construct explicit examples in both situations. In particular, we will show the following:

\begin{thm} Let $M$ be a mean convex $3$-manifold. Let $\Gamma_1$ and $\Gamma_2$ be two simple closed curves in $\partial M$ and let $\Sigma_1$ and $\Sigma_2$ be
two absolutely area minimizing surfaces in $M$ with $\partial \Sigma_i = \Gamma_i$. Let $\Gamma_1\cap\Gamma_2 = \emptyset$, but $\Sigma_1 \cap
\Sigma_2 \neq \emptyset$. Then, either $\Gamma_1$ is not homologous to $\Gamma_2$ in $\partial M$, or $H_2(M)$ is not trivial. Also, there are
examples in both cases.
\end{thm}

Note that the statement of the theorem is equivalent to the result mentioned in the first paragraph \cite{Co}, \cite{Ha}. In particular, in this
setting the statement "$\Sigma_1$ and $\Sigma_2$ are not homologous in $M$ (relative to $\partial M$)" is equivalent to say that either $\partial
\Sigma_1$ is not homologous to $\partial \Sigma_2$, or $\Sigma_1$ and $\Sigma_2$ differ by a homology class of $H_2(M)$. Our examples show that both
situations are possible.

\subsection{Example III-A: $H_2(M)$ is trivial.}

In this part, we will describe an example when the ambient space $M$ has trivial second homology. Let $\Sigma_i$ and $\Gamma_i$ be as in the theorem,
and let $H_2(M)$ be trivial. Then, $\Sigma_1$ and $\Sigma_2$ are homologous in $M$ (rel. $\partial M$) if and only if $\Gamma_1$ and $\Gamma_2$ are
homologous in $\partial M$. This is because if $S$ is a subsurface in $\partial M$ with $\partial S = \Gamma_1 \cup \Gamma_2$, then
$\Sigma_1\cup\Sigma_2 \cup S$ would be a closed surface in $M$. As $H_2(M)$ is trivial, it separates $M$, and hence $\Sigma_1$ and $\Sigma_2$ are
homologous. In the other direction, if $\Sigma_1$ and $\Sigma_2$ are homologous, then they separate a piece $\Omega$ from $M$. Then, $\partial
\Omega$ contains $\Sigma_1$ and $\Sigma_2$. Let $S = \partial \Omega \cap \partial M$. Then, $\partial S = \Gamma_1 \cup \Gamma_2$ which shows that
$\Gamma_1$ and $\Gamma_2$ are homologous in $\partial M$. Hence, we can construct such an example in the trivial homology case only if the disjoint
simple closed curves $\Gamma_1$ and $\Gamma_2$ in $\partial M$ are not homologous in $\partial M$.

Now, we construct the example. Let $B$ be the closed unit ball in $\BR^3$. Let $\gamma_1^+$ be the curve $\partial B \cap \{z=\frac{1}{5}\}$, i.e.
$\gamma_1^+= \{ (x,y,\frac{1}{5})\in \BR^3 \ | \ x^2+y^2 = \frac{24}{25}\}$. Let $\gamma_1^-$ be the curve $\partial B \cap \{z=-\frac{1}{10}\}$.
Similarly, let $\gamma_2^+$ be the curve $\partial B \cap \{z=\frac{1}{10}\}$, and let $\gamma_2^-$ be the curve $\partial B \cap
\{z=-\frac{1}{5}\}$. Since the total area of the disks $D_1^+$ and $D_1^-$ with $\partial D_i^+ =\gamma_i^+$ ($\sim 2\pi$) is greater than the
annulus with boundary $\gamma_1^+\cup\gamma_2^+$ ($\sim \frac{3}{5}\pi$), the absolutely area minimizing surface $A_1$ bounding
$\gamma_1^+\cup\gamma_1^-$ is an annulus, which is a segment of a catenoid. Similarly, let $A_2$ be the absolutely area minimizing surface with
boundary $\gamma_2^+\cup\gamma_2^-$. Clearly, $A_1\cap A_2\neq \emptyset$. Indeed, since $A_2$ is just reflection of $A_1$ with respect to
$xy$-plane, $A_1$ intersects $A_2$ on the $xy$-plane.

Since $\partial A_i=\gamma^+_i\cup\gamma^-_i$ is not connected for each $i$, we cannot use $A_1$ and $A_2$ as counterexamples. Our aim is to add a
bridge $S_i$ to $\partial A_i$ so that new surfaces $\Sigma_i \sim A_i \cup S^i$ has connected boundaries. However, since $\partial B$ is a sphere
and adding a bridge connecting $\gamma_1^+$ and $\gamma_1^-$ without intersecting $\partial A_2$ is impossible. As we want the resulting boundaries
to be disjoint, too, we need to modify the ambient space by adding handles to $B$ to bypass this problem.

\begin{figure}[t]

\relabelbox  {\epsfxsize=5in

\centerline{\epsfbox{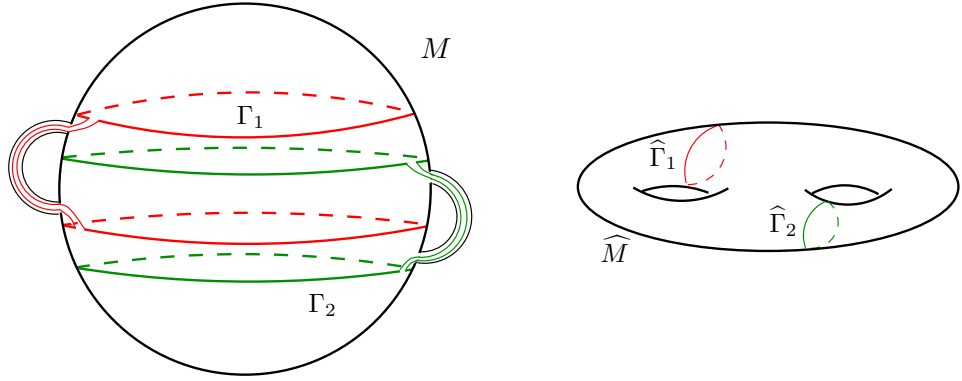}}}

\relabel{1}{\small $\Gamma_1$}

\relabel{2}{\small $\Gamma_2$}

\relabel{3}{$M$}

\relabel{4}{\small $\wh{\Gamma}_1$}

\relabel{5}{\small $\wh{\Gamma}_2$}

\relabel{6}{\small $\wh{M}$}

\endrelabelbox

\caption{\label{intersecting1} \small $M$ is a mean convex $3$-manifold in $\BR^3$. $\Gamma_1$ and $\Gamma_2$ are disjoint simple closed curves in
$\partial M$ (left). Homeomorphic images of $M$, $\Gamma_1$ and $\Gamma_2$ gives the picture in the right, where $\wh{M}$ is a genus $2$ handlebody,
and $\Gamma_1$ and $\Gamma_2$ are non homologous curves in $\partial \wh{M}$.}

\end{figure}

Now, we will add two $1$-handles to $B$. Let $\alpha_1$ be the circular arc with endpoints $x=(0, \sqrt{1-(3/20)^2}, 3/20)$ and
$y=(0,\sqrt{1-(1/20)^2}, -1/20)$, and perpendicular to the the unit sphere $\partial B$ (See Figure \ref{intersecting1} (left)). Let $\alpha_2$ be
the reflection of $\alpha_1$ with respect to the origin. Let $T_1=N_\epsilon(\alpha_1)$, and $T_2=N_\epsilon(\alpha_2)$ be the $1$-handles which we
attach to $B$. Notice that the $3$-manifold $M'=B\cup T_1\cup T_2$ is not mean convex because of the intersections of handles with $\partial B$, i.e.
$\partial T_i\cap \partial B$.

Now, we will modify $M'$ to get a mean convex $3$-manifold. Consider the line $l$ through the point $x$ and perpendicular to the unit sphere
$\partial B$. Parametrize  $l$ by arclength such that $l(0)=x$. Then, let $l(\delta)=x^+$ and $l(-\delta)=x^-$. Let $P^+$ be the plane through $x^+$
and perpendicular to $l$ and similarly, let $P^-$ be the plane through $x^-$ and perpendicular to $l$. Let $\beta^+$ be the round circle $P^+\cap
\partial M'$, and let $\beta^-$ be the round circle $P^-\cap \partial M'$. Now consider the catenoid $C$ with axis $l$ and containing the
circles $\beta^+$ and $\beta^-$. By choosing $\delta$ sufficiently small, we can make sure that the segment $\widehat{C}$ of the catenoid $C$ between
the circles $\beta^+$ and $\beta^-$ lies completely in $M'$. Then, by removing the region between $\widehat{C}$  and $\partial M'$ to $M'$, and do
the same operation at all other $3$ basepoints of the $1$-handles $T_1$ and $T_2$, we get a mean convex manifold $M$. By replacing $\wh{C}$ in the
construction, we can assume $M$ is strictly mean convex, too.

Let $\tau_1$ be a path in $\partial M$ connecting the curves $\gamma_1^+$ and $\gamma_1^-$ through the handle $T_1$. In particular, let $\tau_1$ be
the shorter arc in $\partial M\cap yz$-plane between $\gamma_1^+$ and $\gamma_1^-$, i.e. the endpoints of $\tau_1$ are $(0, \frac{\sqrt{24}}{5},
\frac{1}{5})$ and $(0,\frac{\sqrt{99}}{10}, -\frac{1}{10})$. Similarly, define $\tau_2$ to be the shorter arc in $\partial M\cap yz$-plane between
$\gamma_2^+$ and $\gamma_2^-$ going through the handle $T_2$. Then by using the bridge principle for absolutely area minimizing surfaces (\cite{BC},
Lemma 3.8b), we get the simple closed curve $\Gamma_1$ obtained by putting a thin bridge along $\tau_1$ between $\gamma_1^+$ and $\gamma_1^-$ such
that the absolutely area minimizing surface $\Sigma_1$ in $M$ with $\partial \Sigma_1 = \Gamma_1$ would be the surface obtained from $A_1$ by
attaching a thin strip in $M$ near $\tau_1$, i.e. $\Sigma_1 \sim A_1\cup S_{\tau_1}$ where $S_{\tau_1}$ is a thin strip along $\tau_1$. Similarly, we
get a simple closed curve $\Gamma_2 \sim (\gamma_2^+\cup\gamma_2^-)\sharp \tau_2$ bounding the absolutely area minimizing surface $\Sigma_2$ with
$\Sigma_2 \sim A_2\cup S_{\tau_2}$. By the construction, $\Sigma_1\cap\Sigma_2 \neq \emptyset$. This finishes the example in the trivial homology
case.

Notice that the simple closed curves $\Gamma_1$ and $\Gamma_2$ are not homologous in $\partial M$ as mentioned at the beginning (See Figure
\ref{intersecting1} (right)). The discussion above shows that when $H_2(M)$ is trivial, if $\Gamma_1$ and $\Gamma_2$ are homologous disjoint simple
closed curves in $\partial M$, then the surfaces $\Sigma_1$ and $\Sigma_2$ with $\partial \Sigma_i =\Gamma_i$ must be homologous in $M$. In other
words, if $\Gamma_1\cup \Gamma_2$ separates $\partial M$, then $\Sigma_1\cup\Sigma_2$ must separate $M$. Then, by \cite{Co}, \cite{Ha}, this implies
the absolutely area minimizing surfaces $\Sigma_1$ and $\Sigma_2$ must be disjoint by a simple swaping argument. We rephrase this statement in this
context as follows.

\begin{thm} \cite{Ha, Co} Let $M$ be a strictly mean convex $3$-manifold with trivial $H_2(M)$. Let $\Gamma_1$ and $\Gamma_2$ be two homologous disjoint
simple closed curves in $\partial M$. Then, the absolutely area minimizing surfaces $\Sigma_1$ and $\Sigma_2$ in $M$ with $\partial \Sigma_i =
\Gamma_2$ are disjoint, too.
\end{thm}

\begin{figure}[h]

\relabelbox  {\epsfxsize=3in

\centerline{\epsfbox{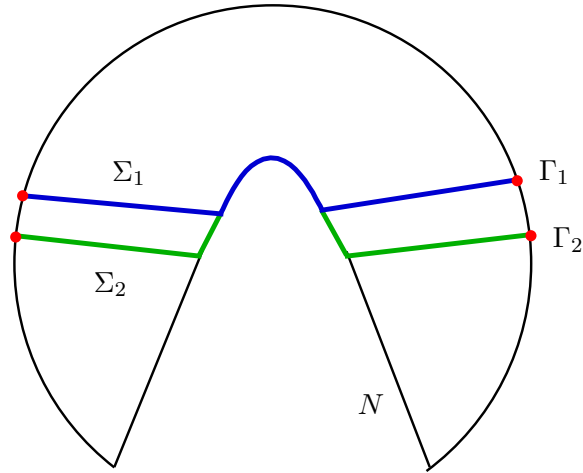}}}

\relabel{1}{$\Gamma_1$}

\relabel{2}{$\Gamma_2$}

\relabel{3}{$\Sigma_1$}

\relabel{4}{$\Sigma_2$}

\relabel{5}{$N$}

\endrelabelbox

\caption{\label{nonmeanconvex} \small In the figure, $N$ represents a non-mean convex $3$-manifold, obtained by removing a large solid cone from a
$3$-ball. $\Gamma_1$ and $\Gamma_2$ are two disjoint simple closed curves in $\partial N$. $\Sigma_1$ and $\Sigma_2$ are the absolutely area
minimizing surfaces in $N$ with $\partial \Sigma_i = \Gamma_i$. Even though $H_2(N)$ is trivial, and $\Gamma_1$ and $\Gamma_2$ are homologous in
$\partial N$, $\Sigma_1 \cap \Sigma_2\neq \emptyset$.}

\end{figure}

\begin{rmk} \label{NMC} In a non-mean convex ambient space $N$, the interaction of the absolutely area minimizing surface and the boundary of the
manifold $\partial N$ can get very complicated. However, for a strictly mean convex $3$-manifold $M$, the absolutely area minimizing surface in $M$
must be away from the boundary $\partial M$ because of the maximum principle. In particular, as it is seen in Figure \ref{nonmeanconvex}, the
absolutely area minimizing surfaces in such a manifold may not be smooth. On the other hand, one can easily construct trivial examples of
intersecting absolutely area minimizing surfaces $\Sigma_1$ and $\Sigma_2$ with disjoint boundaries like in the figure in a non-mean convex
$3$-manifold $N$. To avoid these trivial situations, we construct strictly mean convex examples.
\end{rmk}

\subsection{Example III-B: $H_2(M)$ is nontrivial.}

In this part, we will give another example of intersecting absolutely area minimizing surfaces in $M$ with disjoint boundaries in $\partial M$. In
this example, $H_2(M)$ will be nontrivial, and $\partial M \simeq S^2$. So, any two simple closed curves in $\partial M$ will be homotopic, and
homologous in $\partial M$. Hence, the examples we are going to construct in this part will be very different from the one in previous part in an
essential way.

First we describe the ambient manifold $M$. The ambient manifold $M$ we are going to construct in this part will be very similar to the ambient
manifold in Section 3. Let $T^3$ be the $3$-torus obtained by identifying the opposite faces of the cube with dimensions $[0,1]\times[0,1]\times
[0,1]$. Take the induced flat metric on $T^3$. Let $B^3$ be the cube in the $T^3$ with dimensions $[\delta,1-\delta]\times[\delta,1-\delta]\times
[\delta,1-\delta]$. Define the ambient space $M$ as $T^3-B^3$ where $\delta$ are to be declared later.

\begin{figure}[t]

\relabelbox  {\epsfxsize=5in

\centerline{\epsfbox{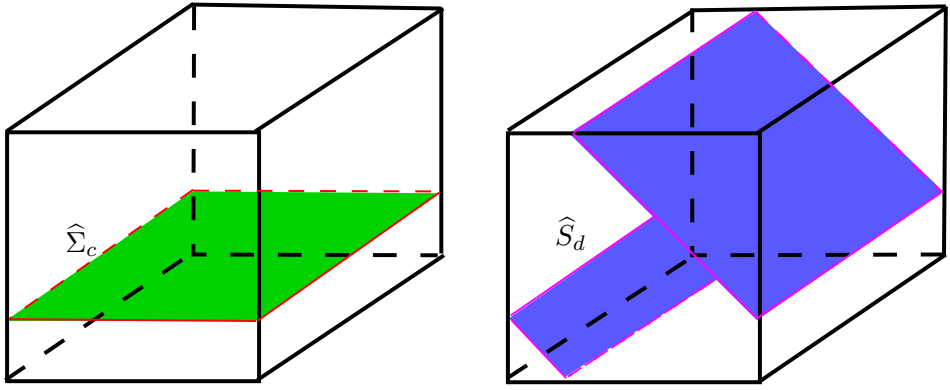}}}

\relabel{1}{$\wh{\Sigma}_c$}

\relabel{2}{$\wh{S}_d$}

\endrelabelbox

\caption{\label{tori} \small {The absolutely area minimizing surfaces $\wh{\Sigma}_c$ and $\wh{S}_d$ in $T^3$.}}

\end{figure}

Let $\wh{\Sigma}_c$ be the be $2$-torus corresponding to the $\{z=c\}$ square in $T^3$ before the identification (See Figure \ref{tori} (left)).
Since the family of minimal tori $\{\wh{\Sigma}_c\}$ foliates $T^3$, $\wh{\Sigma}_c$ is an absolutely area minimizing surface for any $c$. Let
$\Gamma_c$ be the simple closed curve in $\partial M$ with $\Gamma_c = \partial M \cap \wh{\Sigma}_c$ where $\delta<c<1-\delta$. Let $\Sigma_c$ be
the surface obtained by intersecting $\wh{\Sigma}_c$ with $M$ (similar to Figure \ref{horizontal}). Because of the foliation, $\Sigma_c$ is an
absolutely area minimizing surface in $M$ with boundary $\Gamma_c$ unless the competitor disk in $\partial M$ with the same boundary has smaller
area. Hence, an easy computation shows that $\Sigma_c$ is the absolutely area minimizing surface in $M$ with $\partial \Sigma_c = \Gamma_c$ if we
choose $\delta$ such that $\delta < \frac{2-\sqrt2}{4} \sim 0.14$, i.e. $|D|=(1-2\delta)^2 > 1- (1-2\delta)^2 =|\Sigma_c|$ where $D$ is the bottom
(or top) disk in $\partial M$.

Now, we will define the second surface. Let $\wh{S}_d$ be the surface in $T^3$ which is the projection of $y+z=d$-plane in $\BR^3$ (the universal
cover of $T^3$) to $T^3$ via the covering map (See Figure \ref{tori} (right)). In other words, $\wh{S}_d$ is the $2$-torus corresponding to
$\{y+z=d\}$-plane in $T^3$ before identification (in this notation, some piece of $\wh{S}_d$ lives in $\{y+z=d\pm 1\}$-plane).

\begin{figure}[h]
\begin{center}
$\begin{array}{c@{\hspace{.2in}}c}

\relabelbox  {\epsfxsize=2.5in \epsfbox{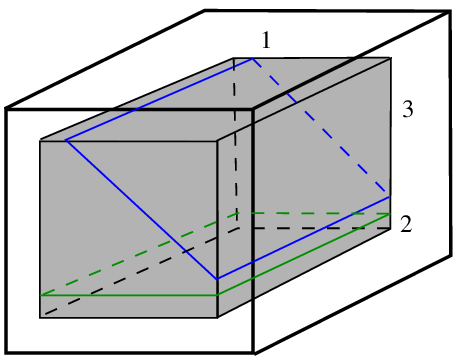}} \relabel{1}{\small $\alpha_d$} \relabel{2}{\small $\Gamma_c$} \relabel{3}{\small $\partial M$}
\endrelabelbox &

\relabelbox  {\epsfxsize=2.5in \epsfbox{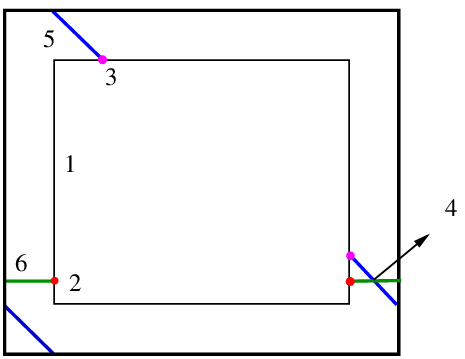}} \relabel{1}{\small $\partial M$} \relabel{2}{\small $\Gamma_c$} \relabel{3}{\small $\alpha_d$}
\relabel{4}{\small $\beta$} \relabel{5}{\small $S_d$} \relabel{6}{\small $\Sigma_c$}  \endrelabelbox \\ [0.4cm]
\end{array}$

\end{center}
\caption{\label{intersecting2} \small In the left, the disjoint simple closed curves $\alpha_d$ (blue curve) and $\Gamma_c$ (green curve) in
$\partial M$ are pictured. In the right, the absolutely area minimizing surfaces $S_d$ and $\Sigma_c$ (right) are intersecting each other along the
simple closed curve $\beta$ even though they have disjoint boundaries $\alpha_d$ and $\Gamma_c$ in $\partial M$ where $M$ is a mean convex
$3$-manifold.}
\end{figure}

Notice that each $\wh{S}_d$ is a minimal surface in $T^3$, as the cover $\{y+z=d\}$-plane is minimal in the universal cover. Since the family of
$2$-tori $\{\wh{S}_d \ | \ d\in [0,1)\}$ foliates $T^3$, each $\wh{S}_d$ is an absolutely area minimizing surface in $T^3$. Let $S_d$ be the surface
which is the intersection of $\wh{S}_d$ and $M$, and $\alpha_d= \wh{S}_d\cap \partial M$ be the simple closed curve(s) in $\partial M$ (See Figure
\ref{intersecting2} (left)). Notice that $\alpha_d$ is one component for $0\leq d <2\delta$ and $1-2\delta<d<1$, and $\alpha_d$ has two components in
$\partial M$ for $2\delta<d<1-2\delta$. From now on, we will assume $0<d<2\delta$. Now, by maximum principle, $S_d$ would be absolutely area
minimizing surface in $M$ with $\partial S_d = \alpha_d$ unless the absolutely area minimizing surface in $M$ with boundary $\alpha_d$ completely
lies in the $\partial M$. Since by our assumption $0\leq d <2\delta$, a simple computation shows that the absolutely area minimizing surface with
boundary $\alpha_d$ cannot lie in $\partial M$. Hence, $S_d$ is an absolutely area minimizing surface in $M$ for $0\leq d <2\delta$.

Now, choose $d=\delta$ and choose $c= \frac{3\delta}{2}$. Then, $\Gamma_c$ and $\alpha_d$ would be disjoint simple closed curves in $\partial M$ (See
Figure \ref{intersecting2} (left)). Moreover, the absolutely area minimizing surfaces $\Sigma_c$ and $S_d$ intersects in a simple closed curve
$\beta$ (See Figure \ref{intersecting2} (right)). Here, $\beta $ corresponds to the line segment between the points $(0, 1- \frac{\delta}{2},
\frac{3\delta}{2})$ and $(1, 1- \frac{\delta}{2}, \frac{3\delta}{2})$ in $T^3$ before identification. This shows the existence of absolutely area
minimizing surfaces $\Sigma_c$ and $S_d$ with nontrivial intersection even though they have disjoint boundaries $\Gamma_c$ and $\alpha_d$ in
$\partial M$.

To get a mean convex example $\wh{M}$, one can follow the steps in Section 4 by modifying the metric on $M$ near $\partial M$. Then suitable
modifications give corresponding absolutely area minimizing surfaces $\wh{\Sigma}_c$ and $\wh{S}_d$ in $\wh{M}$ with the desired properties. This
finishes the second example.

The main difference between the two examples of this section is that in the first example the reason for intersection is not the topological
complexity of the manifold, but the topological difference of the boundaries. In the second example, even though the boundary curves are
topologically "same", the surfaces are in different homological classes which forces the intersection.

\begin{rmk} \label{samebdry} Another interesting question might be "what if the surfaces have the same boundary?". In other words, if $\Sigma_1$ and $\Sigma_2$ are two
absolutely area minimizing surfaces in a mean convex $3$-manifold $M$ with $\partial \Sigma_1 = \partial \Sigma_2 = \Gamma \subset \partial M$. Then,
must $\Sigma_1$ and $\Sigma_2$ be disjoint or not? The answer to this question is clearly "Yes" when $H_2(M)$ is trivial by the discussion at the
beginning of the Section 5.1. However, the answer is "No" when $H_2(M)$ is not trivial.

One can take $\Gamma_c$ and $\Sigma_c$ in the example above. It is possible to construct another absolutely area minimizing surface $T$ with
$\partial T = \Gamma_c$ as follows: Let $T$ be an area minimizing surface in the homology class of $S_d$ (the example above) with $\partial T =
\Gamma_c$. Of course, $T$ is not an absolutely area minimizing surface in $M$ as it is just area minimizing in its homology class. Also,
$T\cap\Sigma_c\neq \emptyset$ by homological reasons. Indeed, the intersection must be in the same homology class with the simple closed curve
$\beta$. Now, we can modify the metric on $M$ near $T$ and away from $\Sigma_c$ so that both $\Sigma_c$ and $T$ are absolutely area minimizing
surface in $M$ with the new metric. In particular, take sufficiently large disk $D$ in $T$ away from the boundary and the intersection, and change
the metric smoothly on a very small neighborhood $N_\epsilon (D)$ of $D$ so that $\Sigma_c$ and $T$ have the same area. Then, we get two intersecting
absolutely area minimizing surfaces with the same boundary.
\end{rmk}

\end{document}